\documentclass[11pt,a4paper]{article}

\usepackage{amsmath,amssymb}
\usepackage{amsmath}
\usepackage{color}
\usepackage[colorlinks,linkcolor=blue,citecolor=red]{hyperref}
\usepackage{tikz}\usepackage{tikz-cd}
\usetikzlibrary{matrix,arrows}

\DeclareMathSymbol{\bbbr}{\mathalpha}{AMSb}{"52}
\DeclareMathSymbol{\bbbc}{\mathalpha}{AMSb}{"52}

\newcommand{\cE}{\mathcal E}

\newcommand\op[1]{\mathop{\rm #1}\nolimits}

\newcommand\com[1]{}
\newcommand\qed{\phantom{\underline{y}}\hfill\hfill$\square$\medskip}

\newtheorem{theorem}{Theorem}

\newtheorem{proposition}[theorem]{Proposition}

\newenvironment{Proof}[1]{\textbf{#1}}{\qed \vspace{5pt}}

\begin{document}

%\title{On integrability of  dispersionless Hirota type equations in 4D}
\title{Second-order PDEs in 4D with half-flat conformal structure}
% in 4D\\ and integrability in higher dimensions}

\author{S. Berjawi$^1$, E.V. Ferapontov$^{1,2}$, B. Kruglikov$^3$, V. Novikov$^1$}
     \date{}
     \maketitle
     \vspace{-5mm}
\begin{center}
$^1$Department of Mathematical Sciences \\ Loughborough University \\
Loughborough, Leicestershire LE11 3TU \\ United Kingdom \\
  \bigskip
  $^2$Institute of Mathematics\\ Ufa Federal Research Centre\\ Russian Academy of Sciences, Russia \\
  \bigskip
$^3$Department of Mathematics and Statistics\\
UiT the Arctic University of Norway\\
Tromso 9037, Norway\\
 [2ex]
e-mails: \\[1ex] \texttt{S.Berjawi@lboro.ac.uk}\\
\texttt{E.V.Ferapontov@lboro.ac.uk}\\
\texttt{Boris.Kruglikov@uit.no} \\
\texttt{V.Novikov@lboro.ac.uk}
\\

\end{center}

\medskip
\begin{abstract}

We study second-order PDEs in 4D for which the conformal structure defined by the characteristic variety
of the equation is half-flat (self-dual or anti-self-dual) on every solution. We prove that this requirement
implies the Monge-Amp\`ere property. Since half-flatness of the conformal structure is equivalent to the
existence of a nontrivial dispersionless Lax pair, our result explains the observation that all known scalar
second-order integrable dispersionless PDEs in dimensions four and higher are of Monge-Amp\`ere type.  
Some partial classification results of Monge-Amp\`ere equations in 4D with half-flat conformal structure 
are also obtained.

\bigskip

\end{abstract}

%\vspace{-7mm}
\newpage

%\tableofcontents

%%%%%%%%%%%%%%%%
\section{Introduction }\label{sec:intro}

Let $u(x)$ be a function of the independent variable $x=(x_1, \dots, x_d)\in M$.
This paper is mainly devoted to the case of dimension four, i.e., $d=\dim M=4$. 
We consider the general class of second-order PDEs  of the form
 \begin{equation}\label{H}
F(x, u, Du, D^2u)=0
 \end{equation}
where $Du=\{u_{\alpha}\}, \ D^2u=\{u_{\alpha \beta}\}$ denote the collection of all first- and 
second-order partial derivatives of $u$ (we use the notation $u_{\alpha}=\partial_{x_{\alpha}}u$, etc). 
 % Here and in what follows the Greek indices take values $\alpha, \beta \in \{1, 2, 3, 4\}$. 
The goal of this paper is to study the implications of integrability on the form of the equation. 
 
%%%%%%%%%%%%%%%%
\subsection{Formulation of the problem and the main results}\label{sec:form}

We will assume that equation (\ref{H}) is {\it non-degenerate} in the sense that its characteristic variety 
defined by the equation
 \begin{equation}\label{quadric}
\sum_{\alpha\leq\beta}\,\frac{\partial F}{\partial u_{\alpha\beta}}\ p_{\alpha} p_{\beta}=0
 \end{equation}
is a non-degenerate quadric on every solution, i.e., for every $j^2u=(x, u, Du, D^2u)$ that satisfies \eqref{H}. 
Similarly, we will say that equation (\ref{H}) has rank $k$ if  quadratic form (\ref{quadric}) has rank $k$ on
every solution. In the non-degenerate case quadratic form (\ref{quadric}) gives rise to the conformal 
structure $[g]$ represented by $g=g_{\alpha \beta} dx_{\alpha} dx_{\beta}$ where $(g_{\alpha \beta})$ 
is the inverse to the matrix of the above quadratic form.
  %$\Bigl(\frac{1+\delta_{\alpha\beta}}2\frac{\partial F}{\partial u_{\alpha \beta}}\Bigr)$
We will be interested in PDEs (\ref{H}) in 4D whose conformal structure is  half-flat 
(self-dual or anti-self-dual) on every solution (see Section~\ref{sec:conf}). 
Examples thereof include heavenly type equations such as the second heavenly equation
 \begin{equation}\label{h2}
u_{13}+u_{24}+u_{11}u_{22}-u_{12}^2=0
 \end{equation}
which governs self-dual Ricci-flat geometry \cite{Plebanski}.
%For the class of dispersionless Hirota type  equations \begin{equation} F(D^2u)=0 \label{H1} \end{equation} this problem was addressed in \cite{FKN} where it was proved that the requirement of half-flatness implies the Monge-Amp\'ere property.
Note that equation \eqref{h2} is of Monge-Amp\`ere type, which is also the case for all other known integrable PDEs of type \eqref{H} for $d\ge4$. 

Recall that a Monge-Amp\`ere equation is a linear combination of all 
possible minors of the Hessian matrix $Hess (u)$ with coefficients being arbitrary functions of the first 
jet $j^1u=(x,u,Du)$. Our main result explains this observation.

 \begin{theorem}\label{T1}
For a non-degenerate PDE (\ref{H}) in dimension $d=4$, half-flatness of the conformal structure $[g]$ 
on every solution implies the Monge-Amp\`ere property. 
 \end{theorem}

It was demonstrated in \cite{FKN} that for a 4D dispersionless Hirota type equation $F(D^2u)=0$, 
half-flatness of $[g]$ is equivalent to  integrability by the method of hydrodynamic reductions, and  implies  the Monge-Amp\`ere property. Thus Theorem~1 is a generalization of the results of \cite{FKN}.
We emphasize that, unlike the method of hydrodynamic reductions, the half-flatness test is fully contact-invariant and applies to general second-order PDEs (\ref{H}).
%The same applies to the existence of a non-trivial dispersionless Lax pair. We can state a somewhat more general result.

It was demonstrated recently in \cite{CalKrug}  that in 4D, half-flatness of $[g]$ is equivalent to the existence of a non-trivial  dispersionless Lax pair in parameter-dependent vector fields. This leads to the following result.

 \begin{theorem}\label{T2}
For a non-degenerate PDE (\ref{H}) in dimension $d=4$, integrability via a non-trivial 
dispersionless Lax pair implies  the Monge-Amp\`ere property. 
 \end{theorem}
 
 Attempts to generalise this result to higher dimensions meet an immediate obstacle: all known multi-dimensional $(d>4)$ PDEs possessing a dispersionless Lax pair are degenerate.  For instance, both the  6-dimensional version
of the second heavenly equation \cite{Takasaki1, Przanovski},
$$
u_{15}+u_{26}+u_{13}u_{24}-u_{14}u_{23}=0,
$$
as well as the 8-dimensional generalisation of the general heavenly equation  \cite{KS},
$$
(u_{x_1y_2} -u_{x_2y_1})(u_{x_3y_4} -u_{x_4y_3}) +(u_{x_2y_3} -u_{x_3y_2})(u_{x_1y_4} -u_{x_4y_1}) +(u_{x_3y_1} u_{x_1y_3})(u_{x_2y_4} -u_{x_4y_2})=0,
$$
have quadratic forms (\ref{quadric})  of rank 4. We however have the following generalisation of Theorem \ref{T2}. 
 
  \begin{theorem}\label{T3}
For a rank 4 PDE (\ref{H}) in any dimension $d\ge4$,  integrability via a non-trivial 
dispersionless Lax pair implies  the Monge-Amp\`ere property. 
 \end{theorem}

It should be noted that in 3D the existence of a nontrivial dispersionless Lax pair does not imply the Monge-Amp\`ere property. 
Generic integrable second-order equations in 3D have transcendental dependence on 2-jets,  see \cite{Fer4}.
 
\medskip

\noindent{\bf Remark 1.} 
The non-degeneracy condition has a simple geometric interpretation in any dimension $d$. 
Without any loss of generality assume that equation (\ref{H}) is of dispersionless Hirota type,  
\vspace{-20pt}

 \begin{equation}\label{H1}
F(D^2u)=0.
 \end{equation}
The corresponding linearised equation is
\vspace{-8pt}

 \begin{equation}\label{H2}
\sum_{\alpha\leq\beta}\frac{\partial F}{\partial u_{\alpha\beta}}\ v_{\alpha \beta}=0,
 \end{equation}
here  $\alpha, \beta \in \{1, \dots, n\}$. Equation (\ref{H1}) can be viewed as the equation  of 
a hypersurface $\cE$ in the Lagrangian Grassmannian $\Lambda$ (locally identified with $n\times n$ 
symmetric matrices $u_{\alpha \beta}$). For any point $L\in\cE$ the tangent space 
$T_L\Lambda\simeq S^2L^*$ (which can also be identified with the space of $n\times n$ 
symmetric matrices $v_{\alpha \beta}$) contains two ingredients:

(a) the Veronese cone $V=\{p\odot p:p\in L^*\}$ of rank 1 matrices; 
 % (flat generalised conformal structure on $\Lambda$). 

(b) the tangent hyperplane $T_L\cE$ defined by equation ({\ref{H2}).  

\noindent In this language, the non-degeneracy of (\ref{H1}) means that $T_L\cE$ is not 
tangential to $V$. %, alternatively,  $T_L\cE$ does not belong to the dual variety $V^*$.
%Indeed, the kernel of the quadratic form \eqref{quadric} on $L^*\simeq T^*M$ is $T_L\cE\cap V$, and so the corank of this form is the dimension of the intersection.

\smallskip

Moreover, consider the bilinear form on $L^*\simeq T^*M$ associated to \eqref{quadric}
at the points of $T_L\cE\cap V$:
 $$
c_F(p,q)=\sum\frac1{1+\delta_{\alpha\beta}}\frac{\partial F}{\partial u_{\alpha\beta}}\,
p_\alpha q_\beta,\quad p,q\in L^*.
 $$
Since $T_pV=\{p\odot q:p,q\in L^*\}$ for $p\neq0$, the kernel of this form is
 $$
\op{Ker} c_F = \{p\in L^*: T_pV\subset T_L\cE\}.
 $$
This is a linear subspace in $L^*$ of dimension equal to $\op{corank}c_F=d-k$.

\medskip

\noindent{\bf Remark 2.}  
Although the converse of Theorem \ref{T1} is not true, our proof implies a somewhat stronger result which says that if equation (\ref{H}) has half-flat conformal structure on every solution then the dependence of 
$F$ on the second-order partial derivatives of $u$ is essentially the same as for 4D integrable symplectic 
Monge-Amp\`ere equations classified in \cite{DF}. That is, freezing the first jet of $u$ (by giving the 
variables $x, u, Du$ arbitrary constant values) we obtain an {\it integrable} symplectic Monge-Amp\`ere 
equation.
Still, a complete classification of PDEs (\ref{H}) with half-flat conformal structure is out of reach at present.

\medskip

In Section \ref{sec:ex} we discuss  examples and partial classification results of PDEs with half-flat 
conformal structure. These are obtained as deformations of some known integrable equations. For example, consider 
the so-called general heavenly equation \cite{Schief},
 $$
\alpha u_{12}u_{34} + \beta u_{13}u_{24}+\gamma u_{23}u_{14} = 0, ~~~ \alpha+\beta+\gamma=0,
 $$
  which,  according to \cite{DF}, is the most generic integrable
symplectic Monge-Amp\`ere equation in 4D.
Note that without the condition  $\alpha+\beta+\gamma=0$ 
the above equation is not integrable, which can be explained on the basis of Theorem \ref{T1}: 
this condition makes it a PDE of Monge-Amp\`ere type. We will show that the general heavenly equation 
possesses the following  deformation with half-flat conformal structure:
 $$
(x_1-x_2)(x_3-x_4)u_{12}u_{34}+(x_3-x_1)(x_2-x_4)u_{13}u_{24}+(x_2-x_3)(x_1-x_4)u_{23}u_{14}=0.
 $$
It will also be shown that these two equations are contact non-equivalent.

\medskip

Proofs of Theorems \ref{T1}, \ref{T2} and \ref{T3} will be given in Section~\ref{Sec3}.
In Section~\ref{Sec4} we discuss  translation non-invariance
of generic integrable second-order PDEs.

%%%%%%%%%%%%%%%%
\subsection{The Monge-Amp\`ere property} \label{sec:MA}

\medskip
 Let us represent equation (\ref{H}) in evolutionary form,
\begin{equation}
% u_{00}=f(u_{0i}, u_{ij}),
u_{11}=f(x, u, Du, u_{1i}, u_{ij})
\label{m0}
\end{equation}
where the Latin indices take values $i, j, k \in \{2, 3, 4\}$.
The proof of Theorem \ref{T1} will require explicit differential constraints for the right-hand side $f$  that are equivalent to the  Monge-Amp\`ere property. These have only been known  in low dimensions \cite{Boillat, Ruggeri, Colin, Gutt}.
In  full generality they were obtained recently in \cite{FKN}. In 4D the Monge-Amp\`ere conditions consist of several groups of equations for $f$.
First of all, for every  $i\in \{2,3,4\}$ one has the relations
 \begin{equation*}\label{E1ab}
f_{u_{ii}}f_{u_{1i}u_{1i}}+f_{u_{ii}u_{ii}}=0, ~~~ f_{u_{1i}}f_{u_{1i}u_{1i}}+2f_{u_{1i}u_{ii}}=0.
 \end{equation*}
Secondly, for every pair of distinct indices $i\ne j\in\{2,3,4\}$ one has the relations
 \begin{gather*}
f_{u_{1j}}f_{u_{1i}u_{1i}}+2f_{u_{1i}}f_{u_{1i}u_{1j}}+2f_{u_{1i}u_{ij}}+2f_{u_{1j}u_{ii}}=0,\label{E2a}\\
f_{u_{ij}}f_{u_{1i}u_{1i}}+2f_{u_{ii}}f_{u_{1i}u_{1j}}+2f_{u_{ii}u_{ij}}=0,\label{E2b}\\
f_{u_{jj}}f_{u_{1i}u_{1i}}+f_{u_{ii}}f_{u_{1j}u_{1j}}+2f_{u_{ij}}f_{u_{1i}u_{1j}}+2f_{u_{ii}u_{jj}}+f_{u_{ij}u_{ij}}=0.\label{E2c}
 \end{gather*}
Furthermore, for every triple of
distinct indices $i\ne j\ne k \in\{2,3,4\}$ one has the relations
 \begin{gather*}
f_{u_{1k}}f_{u_{1i}u_{1j}}+f_{u_{1j}}f_{u_{1i}u_{1k}}+f_{u_{1i}}f_{u_{1j}u_{1k}}+f_{u_{1i}u_{jk}}+f_{u_{1j}u_{ik}}+f_{u_{1k}u_{ij}}=0,\label{E3a}\\
f_{u_{jk}}f_{u_{1i}u_{1i}}+2f_{u_{ik}}f_{u_{1i}u_{1j}}+2f_{u_{ij}}f_{u_{1i}u_{1k}}+2f_{u_{ii}}f_{u_{1j}u_{1k}}+2f_{u_{ii}u_{jk}}+2f_{u_{ij}u_{ik}}=0.
\label{E3b}
 \end{gather*}
Due to the contact invariance of the Monge-Amp\`ere class, this system of 25 relations is invariant
under arbitrary contact transformations.

%%%%%%%%%%%%%%%%
\subsection{Half-flatness and Lax pairs}\label{sec:conf}

In 4D, the key  invariant of a conformal structure  $[g]$ is its Weyl tensor $W$. Let us introduce its 
self-dual and anti-self-dual parts, $W_\pm=\frac{1}{2}(W\pm*W)$, where $*$ is the Hodge star operator
 defined as
$$
*W^i_{jkl}=\frac{1}{2}\sqrt{\det g}\ g^{ia}g^{bc}\epsilon _{ajbd} W^d_{ckl}.
$$
A conformal structure is said to be half-flat (self-dual or anti-self-dual) if $W_-$ or $W_+$ vanishes. Note that the conditions of self-duality and anti-self-duality are  
 % equivalent
switched under the change of orientation.
% \begin{equation*}%\label{self-dual} W=* W, \end{equation*}
 Integrability of the conditions of half-flatness  by the twistor construction is due to Penrose \cite{Penrose}
who observed that half-flatness of $[g]$ is equivalent to the existence of a 3-parameter family of totally null surfaces, see also \cite{DFK} for a more direct treatment of half-flat structures in specially adapted coordinates.
%We also refer to \cite{DFK} for a direct coordinate treatment of conditions of self-duality.
 Thus, for the second heavenly equation (\ref{h2}) we have
%a metric representative is
 $$
g=dx_1dx_3+dx_2dx_4-u_{22}dx_3^2+2u_{12}dx_3dx_4-u_{11}dx_4^2,
 $$
and direct calculation shows that the conformal structure $[g]$ is  half-flat on every solution of (\ref{h2}).

Equations with half-flat conformal structure are known to possess dispersionless Lax pairs, that is,  
vector fields $X, Y$ depending on $j^2u$ and an auxiliary parameter $\lambda$ 
such that the Frobenius integrability condition $[X, Y] \in\text{span}(X,Y)$ holds identically modulo the equation (and its differential consequences). 
 % Taking appropriate linear combinations of $X, Y$ we can always ensure that $[X, Y]=0$ modulo the equation.  
For the second heavenly equation we have
$$
X=\partial_4+u_{22}\partial_2-u_{12}\partial_1+\lambda \partial_1, ~~~
Y=\partial_3-u_{12}\partial_2+u_{22}\partial_1-\lambda \partial_2,
$$
here $\partial_{\alpha}=\frac{\partial}{\partial x_{\alpha}}$.
%, note the absence of terms with $\partial_{\lambda}$.
Integral surfaces of the involutive distribution $\langle X, Y\rangle$
%for various fixed values of the spectral parameter $\lambda$
give totally null surfaces of the corresponding conformal structure $[g]$, thus providing an alternative proof of  its half-flatness.
This is a particular case of the general result of \cite{CalKrug} relating half-flatness 
%of the characteristic conformal structure $[g]$ 
to the existence of a dispersionless Lax pair.
 The main technical
observation of \cite{CalKrug} is that any nontrivial Lax pair must be characteristic, i.e.\ null with respect
to the conformal structure defined by the characteristic variety. 

Lax pairs in vector fields play a key role in the twistor-theoretic technique 
\cite{Atiyah, Penrose, Dun1, Dun2, Strachan1}, dispersionless d-bar method \cite{Bogdanov2} 
and the novel inverse scattering approach \cite{ManSan1, ManSan2} to PDEs 
with half-flat conformal structure.

%%%%%%%%%%%%%%%%
\section{Examples and classification results}\label{sec:ex}

In this section we provide  examples and partial classification results of 
integrable second-order PDEs in 4D. Integrability is manifested via the half-flatness condition. 
 %Further examples we found in \cite{KM} based on the symmetry approach.
All of them are obtained as translationally non-invariant deformations of the following heavenly type Monge-Amp\`ere equations appearing in self-dual Ricci-flat geometry:
\begin{enumerate}
\item $u_{13}+u_{24}+u_{11}u_{22}-u_{12}^2=0$ (second heavenly equation \cite{Plebanski});
\item $u_{13}u_{24}-u_{14}u_{23}=1$ (first heavenly equation \cite{Plebanski});
\item $\alpha u_{12}u_{34} + \beta u_{13}u_{24}+\gamma u_{14}u_{23} = 0$ (general
heavenly equation \cite{Schief}), $\alpha+\beta+\gamma=0$.
\end{enumerate}
In all cases we present dispersionless Lax pairs in $\lambda$-dependent vector fields.
We refer to \cite{KM} for further examples of this kind.

%%%%%%%%%%%%%%%%
\subsection{Equations of the second heavenly type}

 Here we consider equations of the form  $$u_{13}+u_{24}+f(x, u, Du)(u_{11}u_{22}-u_{12}^2)=0$$ 
 where $f$ is  to be determined from the requirement of half-flatness.
The corresponding conformal structure can be represented as
$$
g=dx_1dx_3+dx_2dx_4-fu_{22}dx_3^2+2fu_{12}dx_3dx_4-fu_{11}dx_4^2,
$$
note that $\det g=1$. With the choice $\sqrt {\det g}=1$ one can show that the requirement of vanishing of the anti-self-dual part  $W_-$  leads to a  system of differential constraints for the function $f$ with the general solution
$$
f=c_0(x_1x_3+x_2x_4)+c_1x_1+c_2x_2+c_3x_3+c_4x_4+c_5
$$
where the constants $c_i$ satisfy a single quadratic constraint  $c_1c_3+c_2c_4-c_0c_5=0$ (on the contrary, the vanishing of the self-dual  part $W_+$  leads  to the degenerate case $f=0$).

Modulo % translations and 
linear transformations of the independent variables one can bring $f$ to one of the four normal forms summarised in Table 1 below.
\begin{center}
 \centerline{\footnotesize{Table 1: Equations of the second heavenly type}}
 \medskip
\begin{tabular}{|l|c|}
\hline Function  $f$  & Lax pair $X, Y$ \\
\hline  $f=x_1x_3+x_2x_4$ & $X=\partial_4+(x_1x_3+x_2x_4)(u_{11}\partial_2-u_{12}\partial_1)+\frac{x_1-\lambda x_2}{x_4+\lambda x_3} \partial_1$, \\
&$Y=\partial_3+(x_1x_3+x_2x_4)(u_{22}\partial_1-u_{12}\partial_2)-\frac{x_1-\lambda x_2}{x_4+\lambda x_3} \partial_2$ \\
\hline  $f=x_1$ &  $X=\partial_4+x_1(u_{11}\partial_2-u_{12}\partial_1)+\frac{x_1}{x_4-\lambda} \partial_1$, \\
&$Y=\partial_3+x_1(u_{22}\partial_1-u_{12}\partial_2)-\frac{x_1}{x_4-\lambda} \partial_2$\\
\hline  $f=x_3$ & $X=\partial_4+x_3(u_{11}\partial_2-u_{12}\partial_1)+\frac{\lambda-x_2}{x_3} \partial_1$, \\
&$Y=\partial_3+x_3(u_{22}\partial_1-u_{12}\partial_2)-\frac{\lambda-x_2}{x_3} \partial_2$ \\
\hline $f=1$ & $X=\partial_4+u_{11}\partial_2-u_{12}\partial_1+\lambda \partial_1$, \\
&$Y=\partial_3+u_{22}\partial_1-u_{12}\partial_2-\lambda \partial_2$
 \\
\hline
\end{tabular}
\end{center}
It can be demonstrated that contact symmetry algebras of PDEs corresponding to cases 1 and 2 of Table 1 depend on three arbitrary functions  of two variables. A more detailed analysis revealed that these cases are in fact point-equivalent. The explicit link  is provided by the formulae
$$
X_1=x_1+\frac{x_2x_4}{x_3}, \qquad X_2=\frac{x_2}{x_3}, \qquad X_3=-\frac{1}{x_3}, \qquad X_4=\frac{x_4}{x_3}, \qquad U=x_3 u;
$$
here $u(x_i)$ and $U(X_i)$ satisfy PDEs corresponding to cases 1 and 2, respectively. Similarly, contact symmetry algebras of PDEs correspoding to cases 3 and 4 depend on four arbitrary functions of two variables.  These PDEs  are  related by the point transformation
$$
X_1=x_1x_3, \qquad X_2=\frac{x_2}{x_3^2}, \qquad X_3=-\frac{1}{2x_3^2}, \qquad X_4=x_4, \qquad U=x_3 u+\frac{x_1x_2^2}{2x_3};
$$
here $u(x_i)$ and $U(X_i)$ satisfy PDEs corresponding to cases 3 and 4, respectively. Thus, Table 1 contains  only two   contact non-equivalent cases.

%%%%%%%%%%%%%%%%
\subsection{Equations of the first heavenly type }

Here we consider equations of the form
 $$u_{13}u_{24}-u_{14}u_{23}=f(x, u, Du).$$
The corresponding conformal structure can be represented as
$$
g=u_{13}dx_1dx_3+u_{14}dx_1dx_4+u_{23}dx_2dx_3+u_{24}dx_2dx_4.
$$
Note that by virtue of the equation we have $\det g=\frac{1}{16}f^2$ so that  $\sqrt {\det g}=\frac{1}{4}f$. With this choice of the square root  the requirement of vanishing of the anti-self-dual part  $W_-$  leads to a system of differential constraints for $f$ which we do not present here explicitly due to its complexity (on the contrary, the vanishing of the self-dual  part $W_+$  leads  to the degenerate case $f=0$).
The classification of functions $f$ satisfying the condition $W_-=0$ is performed modulo the following point transformations which leave the class under study form-invariant (this considerably reduces the number of  cases):

(a) changes of variables $x_1\to a(x_1, x_2), \ x_2 \to b(x_1, x_2), \ x_3\to p(x_3, x_4), \ x_4 \to q(x_3, x_4)$;

(b) simultaneous permutations $\{x_1 \leftrightarrow x_3\},\ \{x_2 \leftrightarrow x_4\}$;

(c) translations $u\to u+f(x_1, x_2)+g(x_3, x_4)$;

(d)  rescalings $u\to cu$.
%\begin{center} \begin{tabular}{|l|c|} \hline \centerline{Classification of integrable cases}\\ \hline  $f= \frac{u_1u_3}{(x_2-x_4)^2}$, \quad general case \\  \hline  \small{$f= u_1$} \\  \hline   \small{$f= 1$}, \quad 1st heavenly equation \\ \hline \end{tabular} \end{center}

\noindent The classification results are summarised in  Table 2 below:
\begin{center}
 \centerline{\footnotesize{Table 2: Equations of the first heavenly type}}
 \medskip
\begin{tabular}{|l|c|}
\hline  Function $f$  & Lax pair $X, Y$ \\
\hline $f= \frac{u_1u_3}{(x_2-x_4)^2}$ &   $X=u_{13}\partial_4-u_{14}\partial_3+\frac{u_3(\lambda -x_2)}{(\lambda-x_4)(x_2-x_4)}\partial_1$, \\
& $Y=u_{23}\partial_4-u_{24}\partial_3+\frac{u_3(\lambda -x_2)}{(\lambda-x_4)(x_2-x_4)} \partial_2$ \\
\hline $f= u_1$ &  $X=u_{13}\partial_4-u_{14}\partial_3+(\lambda -x_2)\partial_1$, \\
& $Y=u_{23}\partial_4-u_{24}\partial_3+(\lambda - x_2) \partial_2$  \\
\hline $f=1$ & $X=u_{13}\partial_4-u_{14}\partial_3+\lambda \partial_1$, \\
& $Y=u_{23}\partial_4-u_{24}\partial_3+\lambda \partial_2$
 \\
\hline
\end{tabular}
\end{center}
Contact non-equivalence of the corresponding PDEs can be seen from the structure of their  symmetries: one can show that  contact symmetry algebras of PDEs corresponding to cases 1-3 of Table 2 depend on two, three and four arbitrary functions of two variables, respectively.

%%%%%%%%%%%%%%%%
\subsection{Equations of the general heavenly type   }

Here we classify equations of the form
$$
u_{13}u_{24}-u_{14}u_{23}=f(x, u, Du)(u_{13}u_{24}-u_{12}u_{34}).
$$
In this case the explicit form of $[g]$, as well as the system of differential constraints for $f$ coming from the requirement of half-flatness  are rather lengthy so we skip the details. One can show that the generic such $f$ is given by the formula
$$
f=\frac{(\phi_1-\phi_2)(\phi_3-\phi_4)}{(\phi_1-\phi_4)(\phi_3-\phi_2)}
$$
where $\phi_i=\phi_i(x_i, u_i)$ are arbitrary functions of the indicated variables. If $\phi_i\ne const$ one 
can set  $\phi_i=x_i$ via a suitable contact transformation, thus, normal forms depend on how many functions among $\phi_i$ are non-constant (five different forms altogether). 

PDEs corresponding to such normal forms are contact non-equivalent. Indeed, when $s$ functions 
among $\phi_i$ are constant, the contact symmetry algebra depends on $s$ arbitrary functions 
of two variables, $0\leq s\leq 4$.
The two limiting cases are summarised in   Table 3 below.
%\begin{center}  \centerline{\footnotesize{Table 5: Equations of the general heavenly type}}  \medskip \begin{tabular}{|l|c|} \hline Function  $f$  & Lax pair $X, Y$ \\ \hline $f=\frac{(x_1+x_2)(x_3+x_4)}{(x_1+x_4)(x_3+x_2)}$ &  \\ &  \\ \hline $f=c$ & $X=u_{34}\partial_1-u_{13}\partial_4+\lambda (u_{14}\partial_3-u_{34}\partial_1), ~~Y=u_{23}\partial_4-u_{34}\partial_2+(1-c)\lambda(u_{34}\partial_2-u_{24}\partial_3)$ \\ \hline \end{tabular} \end{center}

%\hline $f=c$ & $X=(1-\lambda) \partial_1-\frac{u_{13}}{u_{34}}\partial_4+\lambda \frac{u_{14}}{u_{34}}\partial_3$,
% $Y=((1-c)\lambda-1)\partial_2+\frac{u_{23}}{u_{34}}\partial_4-(1-c)\lambda \frac{u_{24}}{u_{34}}\partial_3$ \\

\begin{center}
 \centerline{\footnotesize{Table 3: Equations of the general heavenly type}}
 \medskip
\footnotesize{\begin{tabular}{|l|c|}
\hline Function  $f$  & Lax pair $X, Y$ \\
\hline $f=\frac{(x_1-x_2)(x_3-x_4)}{(x_1-x_4)(x_3-x_2)}$ & $X=(x_4-x_3)(x_1+\lambda)u_{34}\partial_1+(x_3-x_1)(x_4+\lambda)u_{13}\partial_4+(x_1-x_4)(x_3+\lambda)u_{14}\partial_3$, \\
 & $Y=(x_4-x_3)(x_2+\lambda)u_{34}\partial_2+(x_3-x_2)(x_4+\lambda)u_{23}\partial_4+(x_2-x_4)(x_3+\lambda)u_{24}\partial_3$\\
%\hline $f=\frac{x_1(x_3+x_4)}{x_3(x_1+x_4)}$ &  \\
\hline $f=c, \ c\ne 0, 1$ & $X=u_{34}\partial_1-u_{13}\partial_4+\lambda (u_{14}\partial_3-u_{34}\partial_1)$, \\
 & $Y=u_{23}\partial_4-u_{34}\partial_2+(1-c)\lambda(u_{34}\partial_2-u_{24}\partial_3)$
\\
\hline
\end{tabular}}
\end{center}
In the last case (which is equivalent to the general heavenly equation) the contact symmetry algebra depends on four arbitrary functions of two variables, while in the first case the contact symmetry algebra involves a certain number of functions of one variable only. This confirms contact non-equivalence of the corresponding equations.

%%%%%%%%%%%%%%%%
\section{Proofs of the main results}\label{sec:proof} \label{Sec3}

 \begin{Proof}{Proof of Theorem \ref{T1}}
is based on  direct computation of the Weyl tensor  of the  conformal structure $[g]$,
\begin{equation}
W_{ijkl}=R_{ijkl}-w_{ik}g_{jl}-w_{jl}g_{ik}+w_{jk}g_{il}+w_{il}g_{jk}=0,
\label{Weyl}
\end{equation}
where $R_{ijkl}= g_{is}R^s_{jkl}$ is the Riemann curvature tensor,
$w_{ij}=\frac12R_{ij}-\frac{R}{12}g_{ij}$ is the 4-dimensional Schouten tensor,
$R_{ij}$ is the Ricci tensor and $R$ is the scalar curvature.  The key calculation (performed in Mathematica) can be split  into several steps:

\begin{itemize}
\item Calculate the Weyl tensor $W$  of the conformal structure $[g]$ defined by the characteristic variety of equation (\ref{m0}). Note that since $[g]$ depends on no more than second-order partial derivatives of  $u$,  components of  $W$ will depend on partial derivatives of $u$ up to the order at most 4. Furthermore, fourth-order partial derivatives of $u$ will enter linearly.

\item Restrict $W$ to a solution of (\ref{m0}), that is, use  (\ref{m0}) and its differential consequences 
to eliminate from $W$ all partial derivatives of $u$ of the form $u_{11}, Du_{11}, D^2u_{11}$.

\item Calculate  the anti-self-duality condition $W_-=0$ (the self-duality  case $W_+=0$ will be essentially the same).

\item Take coefficients of  $W_-$ at the remaining fourth-order partial derivatives of $u$, that is, 
at $u_{ijkl}$, and equate them to zero. This will result in a system of second-order relations for the 
right-hand side $f$ of equation (\ref{m0}), 30 linearly independent equations altogether. 
These relations will involve partial derivatives of $f$ with respect to the variables $u_{ij}$ only.

\item Show that 30 second-order relations for $f$ obtained at the previous step contain all of the 25 relations from Section \ref{sec:MA} that characterise Monge-Amp\`ere equations in 4D. The remaining 5 relations are the additional necessary conditions for integrability, they are equivalent to the requirement that the equation obtained from (\ref{m0}) by freezing the first jet of $u$ (by giving the variables $x, u, Du$ arbitrary constant values) results in an integrable Hirota type equation (we will not use these extra relations in what follows).

\item To simplify the calculations we utilise the fact that the induced action of the contact group
on the space of 1-jets of  $f$  has a unique open orbit (its complement consists of 1-jets of degenerate systems). This property plays a key role in the proof  by allowing one to assume that all sporadic factors depending on first-order derivatives of $f$  that arise in the process of Gaussian elimination at the previous step,  are nonzero. This considerably simplifies the arguments by eliminating non-essential branching.
 % Furthermore, in the verification of various identities involving higher-order partial derivatives of $f$ 
 % one can, without any loss of generality, assign the first-order jet of $f$ any numerical values
 % corresponding to a non-degenerate 1-jet: this often renders otherwise impossible 
 % computations manageable.

More precisely,  calculations are simplified considerably if, after having computed the Weyl tensor $W$ (so that all remaining calculations are entirely algebraic), we use  the  1-jet of  $f$ defined as
$f_{u_{14}}=f_{u_{23}}=1$, while all other  first-order partial derivatives of $f$ are zero (the same 1-jet of $f$ should be substituted into all of the 25 Monge-Amp\`ere relations).
\end{itemize}
%The main steps above follow those of \cite{FKN} for the Hirota type equations, which is possible  because a part of our computation consists of the same symbolic formulae. 
This finishes the proof.
%\noindent The relevant Mathematica programmes are available from the ArXiv version of this paper.
 \end{Proof}

Before passing to the next proof recall that non-triviality of dispersionless Lax pairs in 4D was 
discussed in details in \cite{CalKrug} (in 3D the non-triviality condition of \cite{CalKrug} is somewhat more involved). 
In higher dimensions $d>4$ we  say that a disperionless Lax pair is non-trivial if 
the commutativity of the generators of any modification of the Lax pair which is identical on the equation, holds essentially modulo the equation.

\medskip

\noindent
 \begin{Proof}{Proof of Theorems \ref{T2}, \ref{T3}.} As mentioned in the Introduction, Theorem \ref{T2} follows from the general result of \cite{CalKrug}. As for Theorem \ref{T3}, note that
a generic 4D traveling wave reduction of a multi-dimensional integrable 
PDE of rank 4 will automatically be non-degenerate and integrable, because a generic travelling wave reduction
of a non-trivial disperionless Lax pair is itself non-trivial. Hence the reduced PDE must be of 
 Monge-Amp\`ere type by Theorem \ref{T1} demonstrated above. Indeed, by \cite{CalKrug} 
the existence of a non-trivial Lax pair in 4D yields half-flatness of the conformal structure 
on any solution, so Theorem \ref{T1} applies. Henceforth, since all traveling wave reductions of a multi-dimensional PDE are of Monge-Amp\'ere type, by the argument  similar to that in \cite{FKN}
we conclude that the initial equation in $d$ dimensions should  be of the Monge-Amp\`ere type as well. 
 \end{Proof}

%%%%%%%%%%%%%%%%
\section{Translational invariance for integrable equations}\label{Sec4}

Although some PDEs obtained in Section \ref{sec:ex} contain explicit dependence on the independent 
variables $x=(x_1,\dots,x_d)$, each of them can be put into translationally invariant form by a suitable 
contact transformation. Indeed, the necessary and sufficient condition for this, in general dimension $d$,
is the existence of a $d$-dimensional commutative subalgebra in the algebra of contact symmetries 
of the equation, which acts simply transitively in $J^1M$. This subalgebra is responsible for translations of the 
independent variables and can be contactly mapped to $\langle\partial_{x_1},\dots,\partial_{x_d}\rangle$.

 \begin{proposition}\label{Pro}
There exist nondegenerate integrable PDEs in 3D and 4D that are not contact equivalent to 
any translationally invariant equation.
 \end{proposition}

In other words, there exist examples of integrable PDEs with essential dependence on the independent 
variables. 
%These extend the class of integrable PDEs  considered earlier.
Note that PDEs containing explicit dependence on the independent variables appeared previously,
see e.g.\ \cite{FerKrug,DFK,KM,KP}. However, to our knowledge, the observation that this dependence 
can be essential is apparently new.

\medskip

 \noindent 
 \begin{Proof}{Proof.}
We use the above observation on the existence of $d$-dimensional Abelian subalgebra of
contact symmetries as a necessary condition. We consider the cases $d=3,4$ in turn.

\medskip

\noindent {\bf Case} $\mathbf{d=3.}$ 
The following 3D equation was obtained in \cite{KP} as an  integrable deformation of the Veronese web equation:
\begin{equation}
(x_1-x_2)u_3u_{12}+(x_2-x_3)u_1u_{23}+(x_3-x_1)u_2u_{13}=0.
\label{Ver}
\end{equation}
It has a Lax pair $[X, Y]=0$ where
 $$
X=\partial_1-\frac{x_3-\lambda}{x_1-\lambda}\frac{u_1}{u_3}\partial_3, \qquad 
Y=\partial_2-\frac{x_3-\lambda}{x_2-\lambda}\frac{u_2}{u_3}\partial_3.
 $$
The contact symmetry algebra of equation (\ref{Ver}) is generated by vector fields
 $$
X_f=f(u){\partial}_u, \qquad 
Y_0={\partial_{x_1}}+{\partial _{x_2}}+{\partial _{x_3}}, \qquad 
Y_1=x_1{\partial _{x_1}}+x_2{\partial _{x_2}}+x_3{\partial _{x_3}}.
 $$
One can easily see that this algebra does not contain any three-dimensional commutative subalgebra. 
Indeed, $Y_0$ and $Y_1$ do not commute, though both commute with $X_f$, and in the algebra 
$\langle X_f\rangle$ of vector fields on the line the maximal Abelian subalgebra has dimension 1.
Thus, PDE (\ref{Ver}) is not contact-equivalent to any translationally invariant equation.

\medskip

\noindent {\bf Case} $\mathbf{d=4.}$  
Analogously, in 4D, consider a linear combination of four copies of equation (\ref{Ver}):
 \begin{eqnarray}\nonumber
&&\,\,\,\,a_4[(x_1-x_2)u_3u_{12}+(x_2-x_3)u_1u_{23}+(x_3-x_1)u_2u_{13}]\\ \label{Ver4}
&&+a_3[(x_4-x_2)u_1u_{24}+(x_2-x_1)u_4u_{12}+(x_1-x_4)u_2u_{14}]\\ \nonumber
&&+a_2 [(x_3-x_4)u_1u_{34}+(x_1-x_3)u_4u_{13}+(x_4-x_1)u_3u_{14}]\\ \nonumber
&&+a_1 [(x_3-x_2)u_4u_{23}+(x_2-x_4)u_3u_{24}+(x_4-x_3)u_2u_{34}]=0,
 \end{eqnarray}
where $a_1, \dots, a_4$ are arbitrary constants (each of them can be normalised to
$0$ or $1$ independently).
We emphasize that a linear combination of 3D commuting flows is by no means a 4D integrable PDE 
in general. However, this is the case for {\it linearly degenerate} commuting flows of the Veronese web hierarchy. Indeed, equation (\ref{Ver4}) has a Lax pair $[X, Y]=0$ 
with $X=X_1-b_1X_3$, $Y=X_2-b_2X_3$, where
$$
\begin{array}{c}
X_1=\partial_1-\frac{x_4-\lambda}{x_1-\lambda}\frac{u_1}{u_4}\,\partial_4, \qquad 
X_2=\partial_2-\frac{x_4-\lambda}{x_2-\lambda}\frac{u_2}{u_4}\,\partial_4, \qquad
X_3=\partial_3-\frac{x_4-\lambda}{x_3-\lambda}\frac{u_3}{u_4}\,\partial_4,\\
\ \\
b_1=\frac{x_3-\lambda}{x_1-\lambda} \frac{a_1u_4-a_4u_1}{a_3u_4-a_4u_3}, \qquad
b_2=\frac{x_3-\lambda}{x_2-\lambda} \frac{a_2u_4-a_4u_2}{a_3u_4-a_4u_3}.
\end{array}
$$
The contact symmetry algebra of equation (\ref{Ver4}), when all $a_i\neq0$, 
is generated by the vector fields
 $$
X_f=f(a_1x_1+a_2x_2+a_3x_3+a_4x_4,u){\partial}_u, \qquad 
Y_0=\sum_{i=1}^4 {\partial_{x_i}}, \qquad Y_1=\sum_{i=1}^4 x_i{\partial _{x_i}}.
 $$
Again, $Y_0$ and $Y_1$ do not commute, and both commute with $X_f$ iff $f=f(u)$.
Then, similar to the 3D case, the maximal Abelian subalgebra has dimension 2.
If $f$ contains essential  dependence on the first argument $\sum_{i=1}^4a_ix_i$,
then $X_f$ does not commute with either of $Y_0, Y_1$. Thus if such $f$ is admissible, the maximal
Abelian subalgebra is a part of $\langle X_f\rangle$. This latter Lie algebra contains an
infinite-dimensional commutative subalgebra corresponding to all $f$ depending only on the first 
argument $\sum_{i=1}^4a_ix_i$. However, the prolongation of vector fields $X_f$ to $J^1M$ 
with all admissible $f$ has an orbit $\langle\partial_u,\sum_{i=1}^4a_i\partial_{u_i},
\sum_{i=1}^4u_i\partial_{u_i}\rangle$ of dimension $3$ only.
Thus, PDE (\ref{Ver4}) is not contact-equivalent to any translationally invariant equation.
 \end{Proof}

Let us finally remark that a trick used in the second half of the proof of Proposition \ref{Pro}
does not work in dimensions $d>4$. For instance, in 5D one can form the following linear combination of equations from the Veronese web hierarchy,
 $$
\sum_{(ijklm)=(12345)}\pi_{ijklm}a_{ij}
[(x_k-x_l)u_mu_{kl}+(x_l-x_m)u_ku_{lm}+(x_m-x_k)u_lu_{mk}]=0,
 $$
 where $\pi_{\sigma}$ is the signature of
the permutation $\sigma$. For generic values of $a_{ij}$ this is a non-degenerate determined PDE,
however it is not integrable: a non-trivial dispersionless Lax pair does not exist.
The last claim can be verified directly, but it also follows from the general result of 
\cite[Theorem 1]{CalKrug}: a dispersionless Lax pair of a determined second-order PDE is necessarily characteristic, 
which in our context means that its congruence of 2-planes $\langle X, Y\rangle$ is co-isotropic with respect to 
the canonical conformal structure, and such 2-planes cannot exist in dimensions $d>4$.

%%%%%%%%%%%%%%%%
\section*{Acknowledgements}
We thank  B. Doubrov and M. Pavlov for useful discussions.
This research was  supported by the EPSRC grant  EP/N031369/1, and also 
by the project Pure Mathematics in Norway, funded by the Trond Mohn Foundation 
and the Troms\o\ Research Foundation. 

%%%%%%%%%%%%%%%%

\end{document}